\documentclass[11pt,a4paper]{article}
\usepackage[utf8]{inputenc}
\usepackage{amsmath,amssymb,amsfonts,amsthm}
\usepackage{geometry}
\usepackage{xcolor}
\definecolor{navy}{RGB}{24,43,107}
\definecolor{darkred}{RGB}{160,20,20}
\definecolor{forest}{RGB}{18,118,72}
\definecolor{amberdark}{RGB}{180,95,6}
\definecolor{tealdark}{RGB}{15,118,110}
\definecolor{slateink}{RGB}{71,85,105}
\usepackage{hyperref}
\usepackage{booktabs}
\usepackage{microtype}
\usepackage{tikz}
\usetikzlibrary{arrows.meta,decorations.markings,positioning,calc}
\usepackage{pgfplots}
\pgfplotsset{compat=1.17}

\hypersetup{
    colorlinks=true,
    linkcolor=navy,
    citecolor=darkred,
    urlcolor=navy
}

\newtheorem{theorem}{Theorem}[section]
\newtheorem{lemma}[theorem]{Lemma}
\newtheorem{proposition}[theorem]{Proposition}
\newtheorem{corollary}[theorem]{Corollary}
\theoremstyle{definition}
\newtheorem{definition}[theorem]{Definition}

\theoremstyle{remark}
\newtheorem{remark}[theorem]{Remark}
\newtheorem{problem}[theorem]{Open Problem}

\newcommand{\R}{\mathbb{R}}
\newcommand{\C}{\mathbb{C}}
\newcommand{\HH}{\mathbb{H}}
\newcommand{\OO}{\mathbb{O}}

\newcommand{\F}{\mathbb{F}}
\newcommand{\E}{\mathbb{E}}
\newcommand{\Prob}{\mathbb{P}}
\newcommand{\calA}{\mathcal{A}}
\newcommand{\calB}{\mathcal{B}}
\newcommand{\calT}{\mathcal{T}}
\newcommand{\ImO}{\operatorname{Im}}
\newcommand{\ReO}{\operatorname{Re}}
\newcommand{\Aut}{\operatorname{Aut}}
\newcommand{\GL}{\operatorname{GL}}

\newcommand{\Span}{\operatorname{span}}

\title{\textbf{Catalan-Averaged Products in the Octonions:\\[0.2em]
\Large Quaternionic Rigidity, Exact Bilinear Recurrence,\\
and Non-Associative Exponential Decay}}

\author{
\textbf{Terence R. Smith}\thanks{Email: \texttt{smithtr@cs.ucsb.edu}}\\
\small Department of Computer Science,\\
\small University of California, Santa Barbara, CA 93106, USA
}

\date{November 2026}

\begin{document}

\maketitle

\begin{abstract}
Because the real octonions $\OO$ form a non-associative normed division algebra, the product of $N$ octonions $x_1, \dots, x_N \in \OO$ depends on the choice of full binary parenthesization tree $T \in \calB_N$, of which there are $C_{N-1} = \frac{1}{N}\binom{2N-2}{N-1}$ Catalan configurations. While norm multiplicativity forces every individual parenthesized product of unit octonions to retain unit norm, their uniform Catalan average $\calA_N(x_1, \dots, x_N) := C_{N-1}^{-1}\sum_{T \in \calB_N} \mathrm{prod}_T(x_1, \dots, x_N)$ experiences destructive interference governed purely by non-associativity. We establish four main results. First, by bilinearity of octonionic multiplication, $\calA_N$ satisfies an exact $O(N^3)$ dynamic-programming recurrence with Catalan split weights, enabling exact evaluation over $C_{N-1} > 10^{44}$ trees at $N = 80$. Second, we prove a Universal Quaternionic Rigidity Theorem: for all $N \ge 1$ and all nonzero $x_1, \dots, x_N \in \OO$, $\|\calA_N(x_1, \dots, x_N)\| = \prod_{i=1}^N \|x_i\|$ if and only if $\{x_1, \dots, x_N\}$ lies inside a single associative quaternionic subalgebra $\HH \subset \OO$; for Fano basis words, non-associative level sets are quantized in multiples of $168 = |\GL(3,\F_2)|$, and exact cancellation $\calA_N = 0$ is impossible whenever $N = 2^k$. Third, we prove an exact isotropic $N = 3$ moment identity tying $\E[\|\calA_3\|^2]$ to the octonionic associator Casimir, yielding $25/49$ on $S^6$ and Fano words, $5/9$ on the transverse $6$-plane $e_7^\perp$, and $43/64$ on $S^7$. Fourth, combining Minkowski's inequality, norm multiplicativity, subword independence, and singularity analysis of a defect-perturbed Catalan generating function, we prove rigorously that $\E[\|\calA_N\|^2]$ decays exponentially and $\|\calA_N\| \to 0$ exponentially almost surely for i.i.d.\ random unit octonions. High-precision $O(N^3)$ computations through $N = 80$ reveal sharp quenched non-associative Lyapunov exponents $\lambda \approx 0.344$ (Fano basis), $0.268$ ($S^6$), $0.244$ ($e_7^\perp$), and $0.222$ ($S^7$), compared to $\lambda = 0$ on $\HH$.
\end{abstract}

\vspace{0.6em}

\bigskip

%==============================================================================
\section{Introduction}
\label{sec:intro}
%==============================================================================

In the classical theory of random matrix products initiated by Furstenberg and Kesten~\cite{furstenberg1960} and Furstenberg~\cite{furstenberg1963}, exponential growth or decay rates of $\|X_N \cdots X_1\|$ are governed by non-commutativity and linear operator stretching. A qualitatively different phenomenon arises at the top of Hurwitz's tower of real normed division algebras $\R \subset \C \subset \HH \subset \OO$~\cite{hurwitz1898, baez2002}. Because the octonions $\OO$ satisfy the Euclidean composition identity $\|xy\| = \|x\|\,\|y\|$, every individual product of $N$ unit octonions $x_1, \dots, x_N \in S^7 \subset \OO$ has unit norm regardless of how factors are grouped. Unlike the associative algebras $\R, \C, \HH$, however, an $N$-fold octonionic word $w = (x_1, \dots, x_N)$ admits
\begin{equation}
\label{eq:catalan_def}
C_{N-1} = \frac{1}{N}\binom{2N-2}{N-1} \sim \frac{4^{N-1}}{\sqrt{\pi}\,N^{3/2}}
\end{equation}
distinct full binary parenthesization trees $T \in \calB_N$~\cite{stanley2015}, each yielding a unit octonion $\mathrm{prod}_T(x_1, \dots, x_N) \in S^7$.

This motivates studying the uniform \emph{Catalan-averaged product}
\begin{equation}
\label{eq:cat_avg_def}
\calA_N(x_1, \dots, x_N) := \frac{1}{C_{N-1}} \sum_{T \in \calB_N} \mathrm{prod}_T(x_1, \dots, x_N) \in \OO.
\end{equation}
On $\R, \C$, and $\HH$, associativity forces all $C_{N-1}$ summands to coincide, so $\|\calA_N(x_1, \dots, x_N)\| = \prod_{i=1}^N \|x_i\|$ identically. On $\OO$, by contrast, different parenthesizations interfere destructively whenever the factors fail to generate an associative subalgebra. The purpose of this paper is to analyze the algebraic, combinatorial, and probabilistic behavior of $\calA_N$ rigorously. (The problem was originally motivated by an open question on associative phase dispersion in discrete octonionic sequence models~\cite{smith2026octonion}, but the analysis presented here is entirely self-contained pure mathematics.)

Our main contributions are as follows:
\begin{enumerate}
    \item \textbf{Exact $O(N^3)$ Bilinear Evaluation (Theorem~\ref{thm:bilinear_dp}):} Although $|\calB_N| = C_{N-1}$ grows exponentially as $\Theta(4^N N^{-3/2})$, bilinearity of multiplication in any non-associative algebra implies that the subword averages $\calA(i, j)$ satisfy a closed $O(N^3)$ dynamic-programming recurrence weighted by the Catalan split distribution $w_{L,m} = C_{m-1}C_{L-m-1}/C_{L-1}$. This enables exact evaluation of $\calA_N$ at $N = 80$ (averaging over $C_{79} \approx 2.89 \times 10^{44}$ trees) in milliseconds without overflow.
    \item \textbf{Universal Quaternionic Rigidity (Theorem~\ref{thm:rigidity}):} For every $N \ge 1$ and arbitrary nonzero $x_1, \dots, x_N \in \OO$, we prove that $\|\calA_N(x_1, \dots, x_N)\| = \prod_{i=1}^N \|x_i\|$ \emph{if and only if} $\{x_1, \dots, x_N\}$ is contained in a single associative quaternionic subalgebra $\HH \subset \OO$. For words in the Fano basis $\{e_1, \dots, e_7\}$, we show further that all non-quaternionic level sets have cardinalities divisible by $168 = |\GL(3,\F_2)|$ (Proposition~\ref{prop:fano_orbit}) and prove via Kummer's theorem that exact zero cancellation $\calA_N(w) = 0$ is impossible whenever $N = 2^k$ (Corollary~\ref{cor:parity}).
    \item \textbf{Exact $N = 3$ Associator Casimir Identity (Theorem~\ref{thm:n3_exact}):} For three i.i.d.\ random unit octonions with isotropic covariance $(1/d)I_d$, we prove $\E[\|\calA_3\|^2] = 1 - \frac{1}{4d^3}\sum_{a,b,c}\|[e_a,e_b,e_c]\|^2$, yielding the exact rational values $25/49$ on $S^6$ and uniform Fano words, $5/9$ on the transverse $6$-plane $e_7^\perp$, and $43/64$ on $S^7$.
    \item \textbf{Rigorous Exponential Decay and Almost-Sure Convergence (Theorem~\ref{thm:exp_decay}):} Applying Minkowski's inequality, norm multiplicativity, and independence of disjoint subwords yields a nonlinear Catalan majorant recurrence $a_N \le \sum_{k=1}^{N-1} w_{N,k}\,a_k\,a_{N-k}$ for $a_N := (\E[\|\calA_N\|^2])^{1/2}$ with $a_3 < 1$. Analyzing its defect-perturbed Catalan generating function $B^*(z) = \frac{1}{2}(1 - \sqrt{P_K(z)})$ via Pringsheim's theorem and Flajolet--Sedgewick singularity analysis~\cite{flajolet2009} proves rigorously that $\E[\|\calA_N\|^2]$ decays exponentially and that $\|\calA_N(x_1, \dots, x_N)\| \to 0$ exponentially almost surely, with strictly positive lower bounds on the non-associative Lyapunov exponent.
    \item \textbf{High-Precision Numerical Lyapunov Spectrum (Section~\ref{sec:numerics}):} Exact $O(N^3)$ computations through $N = 80$ confirm clean exponential decay across all four ensembles, giving quenched exponents $\lambda_{\mathrm{Fano}} \approx 0.344$, $\lambda_{S^6} \approx 0.268$, $\lambda_{V_\perp} \approx 0.244$, and $\lambda_{S^7} \approx 0.222$, strictly ordered by associator strength and exhibiting a small annealed--quenched gap.
\end{enumerate}

The paper is organized as follows. Section~\ref{sec:prelim} reviews the octonion algebra, Fano plane conventions, and Catalan trees. Section~\ref{sec:dp} establishes the $O(N^3)$ bilinear recurrence and sub-unitary contraction property. Section~\ref{sec:rigidity} proves the Universal Quaternionic Rigidity Theorem and presents exact exhaustive Fano-word enumerations through $N = 6$. Section~\ref{sec:exp_decay} proves the $N = 3$ associator identity and the generating-function exponential decay theorem. Section~\ref{sec:numerics} reports numerical measurements through $N = 80$, and Section~\ref{sec:open} outlines open mathematical problems.

%==============================================================================
\section{Preliminaries}
\label{sec:prelim}
%==============================================================================

\subsection{The Octonions and the Fano Plane}

Let $\OO \cong \R^8$ denote the real octonions with canonical orthonormal basis $\{e_0, e_1, \dots, e_7\}$, where $e_0 = 1$ is the multiplicative identity and $e_a^2 = -1$ for $a \in \{1, \dots, 7\}$. Writing $\OO = \R \cdot 1 \oplus \ImO(\OO)$ with $\ImO(\OO) = \Span_\R\{e_1, \dots, e_7\}$, the product of imaginary basis elements is given by
\begin{equation}
\label{eq:oct_mult}
e_a e_b = -\delta_{ab} + \sum_{c=1}^7 \epsilon_{abc}\,e_c \qquad (a, b \in \{1, \dots, 7\}),
\end{equation}
where $\epsilon_{abc}$ is totally antisymmetric and $\epsilon_{abc} = +1$ along the seven cyclically oriented lines of the Fano plane $\mathbb{PG}(2, \F_2)$:
\begin{equation}
\label{eq:fano_triads}
\calT = \big\{(1,2,4),\; (2,3,5),\; (3,4,6),\; (4,5,7),\; (5,6,1),\; (6,7,2),\; (7,1,3)\big\}.
\end{equation}
The Euclidean inner product is $\langle x, y \rangle = \sum_{\mu=0}^7 x_\mu y_\mu = \ReO(x \bar{y})$, where conjugation maps $1 \mapsto 1$ and $e_a \mapsto -e_a$ ($a \ge 1$). The induced norm $\|x\| = \sqrt{x\bar{x}}$ is strictly multiplicative:
\begin{equation}
\label{eq:norm_mult}
\|xy\| = \|x\|\,\|y\| \qquad \text{for all } x, y \in \OO.
\end{equation}
In particular, $\OO$ has no zero divisors, and left and right multiplication operators $L_a(x) := ax$, $R_a(x) := xa$ are invertible linear isometries (up to scaling by $\|a\|$) for every $a \in \OO \setminus \{0\}$.

The \emph{associator} $[x, y, z] := (xy)z - x(yz)$ is trilinear and totally antisymmetric (alternating) in $(x, y, z)$~\cite{schafer1966, baez2002}. By Artin's theorem~\cite{zorn1930, schafer1966}, any two elements $x, y \in \OO$ generate an associative subalgebra of $\OO$ (isomorphic to $\R$, $\C$, or $\HH$), and $\OO$ satisfies the Moufang identities~\cite{moufang1935}. Every quaternionic subalgebra $\HH_0 \subset \OO$ yields a Cayley--Dickson orthogonal decomposition $\OO = \HH_0 \oplus \HH_0 \ell$ for any unit vector $\ell \in \HH_0^\perp$, with multiplication
\begin{equation}
\label{eq:cayley_dickson}
(a + b\ell)(c + d\ell) = (ac - \bar{d}b) + (da + b\bar{c})\ell \qquad (a, b, c, d \in \HH_0).
\end{equation}
On the canonical basis, the seven Fano lines in~\eqref{eq:fano_triads} correspond to seven coordinate quaternionic subalgebras $\HH_{abc} = \Span_\R\{1, e_a, e_b, e_c\}$ (Figure~\ref{fig:fano_dp}(a)). Two distinct units $e_a, e_b$ lie on a unique Fano line $(a,b,c) \in \calT$; if $d \in \{a,b,c\}$ then $[e_a, e_b, e_d] = 0$, whereas if $d \notin \{a,b,c\}$ then $(e_a e_b)e_d = -e_a(e_b e_d)$ and $\|[e_a, e_b, e_d]\|^2 = 4$.

\subsection{Full Binary Trees and Catalan-Averaged Products}

For $N \ge 1$, let $\calB_N$ denote the set of full (rooted, ordered) binary trees with $N$ leaves, so $|\calB_N| = C_{N-1}$ with $C_0 = 1, C_1 = 1, C_2 = 2, C_3 = 5, C_4 = 14, C_5 = 42, \dots$~\cite{stanley2015}. Every tree $T \in \calB_N$ ($N \ge 2$) decomposes uniquely at its root into an ordered pair $(T_L, T_R) \in \calB_k \times \calB_{N-k}$ for a unique split index $k \in \{1, \dots, N-1\}$.

\begin{definition}[Catalan-Averaged Product]
Given a finite-dimensional real or complex algebra $\mathbb{A}$ and an $N$-tuple $w = (x_1, \dots, x_N) \in \mathbb{A}^N$, each $T \in \calB_N$ defines an evaluated product $\mathrm{prod}_T(w) \in \mathbb{A}$ recursively by $\mathrm{prod}_{\bullet}(x_1) := x_1$ and
\begin{equation}
\mathrm{prod}_{(T_L, T_R)}(x_1, \dots, x_N) := \mathrm{prod}_{T_L}(x_1, \dots, x_k)\,\mathrm{prod}_{T_R}(x_{k+1}, \dots, x_N).
\end{equation}
The \emph{Catalan-averaged product} of $w$ is
\begin{equation}
\label{eq:cat_avg}
\calA_N(x_1, \dots, x_N) := \frac{1}{C_{N-1}} \sum_{T \in \calB_N} \mathrm{prod}_T(x_1, \dots, x_N).
\end{equation}
\end{definition}

For $N = 1, 2$ we have $\calA_1(x_1) = x_1$ and $\calA_2(x_1, x_2) = x_1 x_2$; for $N = 3$,
\begin{equation}
\label{eq:A3_def}
\calA_3(x_1, x_2, x_3) = \frac{1}{2}\big((x_1 x_2) x_3 + x_1 (x_2 x_3)\big).
\end{equation}

%==============================================================================
\section{Bilinear Catalan Recurrence and Exact $O(N^3)$ Evaluation}
\label{sec:dp}
%==============================================================================

Naive evaluation of~\eqref{eq:cat_avg} sums over $C_{N-1} = \Theta(4^N N^{-3/2})$ trees, becoming intractable past $N \approx 16$. The following theorem reduces exact evaluation to cubic time in any non-associative algebra.

\begin{theorem}[Bilinear Catalan Recurrence]
\label{thm:bilinear_dp}
Let $\mathbb{A}$ be any (not necessarily associative) algebra over $\R$ or $\C$, and let $w = (x_1, \dots, x_N) \in \mathbb{A}^N$. For $0 \le i < j \le N$, denote the contiguous subword average by
\begin{equation}
\calA(i, j) := \calA_{j-i}(x_{i+1}, \dots, x_j).
\end{equation}
Then $\calA(i, i+1) = x_{i+1}$, and for every subword length $L = j - i \ge 2$,
\begin{equation}
\label{eq:dp_recurrence}
\calA(i, j) = \sum_{k=i+1}^{j-1} w_{L,\, k-i}\,\calA(i, k)\,\calA(k, j),
\end{equation}
where the Catalan split weights
\begin{equation}
\label{eq:split_weights}
w_{L, m} := \frac{C_{m-1}\,C_{L-m-1}}{C_{L-1}} \qquad (1 \le m \le L-1)
\end{equation}
form a probability distribution on $\{1, \dots, L-1\}$ ($\sum_{m=1}^{L-1} w_{L,m} = 1$, $w_{L,m} > 0$). Consequently, $\calA_N(x_1, \dots, x_N) = \calA(0, N)$ can be computed exactly in $\frac{1}{6}N(N^2 - 1) = O(N^3)$ algebra multiplications and $O(N^2)$ storage.
\end{theorem}

\begin{proof}
Let $S(i, j) := C_{j-i-1}\,\calA(i, j) = \sum_{T \in \calB_{j-i}} \mathrm{prod}_T(x_{i+1}, \dots, x_j)$. Partitioning $\calB_{j-i}$ according to the root split $k \in \{i+1, \dots, j-1\}$ gives a bijection $\calB_{j-i} \cong \bigsqcup_{k=i+1}^{j-1} (\calB_{k-i} \times \calB_{j-k})$. By bilinearity (distributivity) of multiplication in $\mathbb{A}$,
\begin{align}
S(i, j) &= \sum_{k=i+1}^{j-1} \sum_{T_L \in \calB_{k-i}} \sum_{T_R \in \calB_{j-k}} \mathrm{prod}_{T_L}(x_{i+1}, \dots, x_k)\,\mathrm{prod}_{T_R}(x_{k+1}, \dots, x_j) \nonumber\\
&= \sum_{k=i+1}^{j-1} \left(\sum_{T_L \in \calB_{k-i}} \mathrm{prod}_{T_L}(x_{i+1}, \dots, x_k)\right)\!\left(\sum_{T_R \in \calB_{j-k}} \mathrm{prod}_{T_R}(x_{k+1}, \dots, x_j)\right) \nonumber\\
&= \sum_{k=i+1}^{j-1} S(i, k)\,S(k, j).
\end{align}
Dividing by $C_{L-1}$ and substituting $S(i, k) = C_{k-i-1}\calA(i, k)$ and $S(k, j) = C_{j-k-1}\calA(k, j)$ yields~\eqref{eq:dp_recurrence}. The identity $\sum_{m=1}^{L-1} w_{L,m} = 1$ is Segner's Catalan recurrence $C_{L-1} = \sum_{m=1}^{L-1} C_{m-1} C_{L-m-1}$. Evaluating $\calA(i, j)$ in increasing order of length $L = 2, \dots, N$ requires $\sum_{L=2}^N (N - L + 1)(L - 1) = \binom{N+1}{3} = \frac{1}{6}N(N^2-1)$ algebra products.
\end{proof}

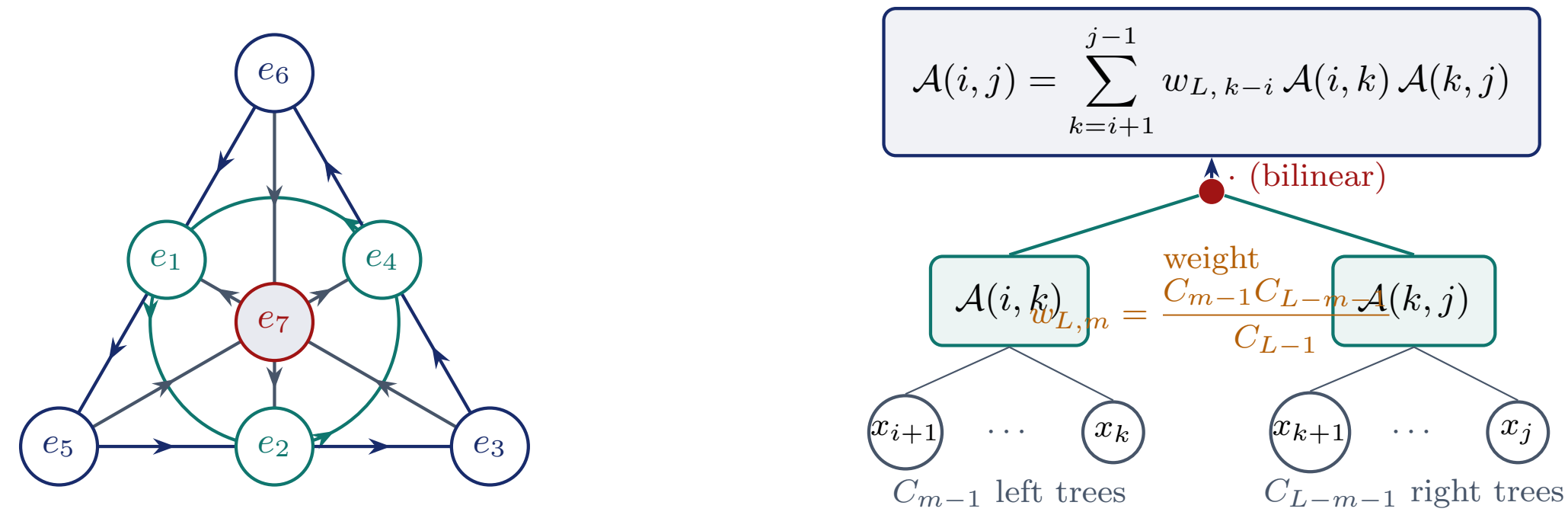
\begin{remark}[Numerical Stability]
Working directly with the normalized recurrence~\eqref{eq:dp_recurrence} rather than unnormalized sums $S(i, j)$ completely avoids floating-point overflow of $C_{N-1} \sim 4^{N-1}/(\sqrt{\pi}N^{3/2})$, since every weight $w_{L,m}$ lies in $(0, 1]$ and every intermediate entry satisfies $\|\calA(i, j)\| \le 1$ for unit inputs. The geometric structure of the Fano plane and the bilinear Catalan root-split factorization are illustrated in Figure~\ref{fig:fano_dp}.
\end{remark}

\begin{figure}[htbp]
\centering
\begin{minipage}[c]{0.47\textwidth}
\centering
\begin{tikzpicture}[
    scale=1.08,
    midarrow/.style={
        postaction={decorate},
        decoration={markings, mark=at position 0.54 with {\arrow{Stealth[length=5.5pt,width=4pt]}}}
    },
    fnode/.style={
        circle, draw=navy, fill=white, thick,
        minimum size=19pt, inner sep=0pt, font=\footnotesize\bfseries, text=navy
    }
]
    % Equilateral triangle vertices and midpoints (radius R = 2.0)
    \coordinate (A) at (90:2.0);      % e_6
    \coordinate (B) at (210:2.0);     % e_5
    \coordinate (C) at (330:2.0);     % e_3
    \coordinate (MAB) at (150:1.0);   % e_1
    \coordinate (MBC) at (270:1.0);   % e_2
    \coordinate (MCA) at (30:1.0);    % e_4
    \coordinate (O) at (0,0);         % e_7

    % Inscribed circle: (1,2,4) counterclockwise
    \draw[thick, tealdark,
        postaction={decorate},
        decoration={markings,
            mark=at position 0.16 with {\arrow{Stealth[length=5.5pt,width=4pt]}},
            mark=at position 0.50 with {\arrow{Stealth[length=5.5pt,width=4pt]}},
            mark=at position 0.83 with {\arrow{Stealth[length=5.5pt,width=4pt]}}
        }] (0,0) circle (1.0);

    % Outer sides: counterclockwise 6->1->5, 5->2->3, 3->4->6
    \draw[thick, navy, midarrow] (A) -- (MAB);
    \draw[thick, navy, midarrow] (MAB) -- (B);
    \draw[thick, navy, midarrow] (B) -- (MBC);
    \draw[thick, navy, midarrow] (MBC) -- (C);
    \draw[thick, navy, midarrow] (C) -- (MCA);
    \draw[thick, navy, midarrow] (MCA) -- (A);

    % Medians inward through center e_7: 6->7->2, 5->7->4, 3->7->1
    \draw[thick, slateink, midarrow] (A) -- (O);
    \draw[thick, slateink, midarrow] (O) -- (MBC);
    \draw[thick, slateink, midarrow] (B) -- (O);
    \draw[thick, slateink, midarrow] (O) -- (MCA);
    \draw[thick, slateink, midarrow] (C) -- (O);
    \draw[thick, slateink, midarrow] (O) -- (MAB);

    % Nodes
    \node[fnode] at (A) {$e_6$};
    \node[fnode] at (B) {$e_5$};
    \node[fnode] at (C) {$e_3$};
    \node[fnode, draw=tealdark, text=tealdark] at (MAB) {$e_1$};
    \node[fnode, draw=tealdark, text=tealdark] at (MBC) {$e_2$};
    \node[fnode, draw=tealdark, text=tealdark] at (MCA) {$e_4$};
    \node[fnode, fill=navy!10, draw=darkred, text=darkred] at (O) {$e_7$};

    \node[font=\footnotesize\itshape, text=slateink] at (0,-1.72) {(a) Oriented Fano plane $\mathbb{PG}(2,\F_2)$};
\end{tikzpicture}
\end{minipage}
\hfill
\begin{minipage}[c]{0.50\textwidth}
\centering
\begin{tikzpicture}[
    scale=0.95,
    box/.style={
        rectangle, rounded corners=3pt, draw=navy, fill=navy!6,
        thick, minimum height=22pt, inner xsep=7pt, font=\footnotesize
    },
    subbox/.style={
        rectangle, rounded corners=3pt, draw=tealdark, fill=tealdark!8,
        thick, minimum height=22pt, inner xsep=6pt, font=\footnotesize
    },
    leaf/.style={
        circle, draw=slateink, fill=white, thick,
        minimum size=15pt, inner sep=0pt, font=\scriptsize
    }
]
    % Root average node
    \node[box, minimum width=128pt] (root) at (0, 2.35) {$\calA(i, j) = \displaystyle\sum_{k=i+1}^{j-1} w_{L,\,k-i}\,\calA(i, k)\,\calA(k, j)$};

    % Multiplication dot
    \node[circle, fill=darkred, inner sep=2.2pt, label={[right, font=\scriptsize, text=darkred]$\cdot$ (bilinear)}] (mul) at (0, 1.35) {};
    \draw[thick, navy, -{Stealth[length=5pt]}] (mul) -- (root);

    % Left and right subword averages
    \node[subbox] (left) at (-1.85, 0.35) {$\calA(i, k)$};
    \node[subbox] (right) at (1.85, 0.35) {$\calA(k, j)$};

    \draw[thick, tealdark] (left.north) -- (mul);
    \draw[thick, tealdark] (right.north) -- (mul);

    % Weight annotation
    \node[font=\scriptsize, text=amberdark, align=center] at (0, 0.35) {weight\\$w_{L,m} = \dfrac{C_{m-1}C_{L-m-1}}{C_{L-1}}$};

    % Subtrees to leaves
    \node[leaf] (x1) at (-2.8, -0.85) {$x_{i+1}$};
    \node[font=\scriptsize, text=slateink] at (-1.85, -0.85) {$\cdots$};
    \node[leaf] (xk) at (-0.9, -0.85) {$x_k$};

    \node[leaf] (xk1) at (0.9, -0.85) {$x_{k+1}$};
    \node[font=\scriptsize, text=slateink] at (1.85, -0.85) {$\cdots$};
    \node[leaf] (xj) at (2.8, -0.85) {$x_j$};

    \draw[slateink, thin] (left.south) -- (x1.north);
    \draw[slateink, thin] (left.south) -- (xk.north);
    \draw[slateink, thin] (right.south) -- (xk1.north);
    \draw[slateink, thin] (right.south) -- (xj.north);

    % Braces / complexity summary
    \node[font=\scriptsize, text=slateink, align=center] at (-1.85, -1.42) {$C_{m-1}$ left trees};
    \node[font=\scriptsize, text=slateink, align=center] at (1.85, -1.42) {$C_{L-m-1}$ right trees};
    \node[font=\footnotesize\itshape, text=slateink] at (0, -2.02) {(b) $O(N^3)$ Catalan split factorization};
\end{tikzpicture}
\end{minipage}
\caption{(a)~The oriented Fano plane corresponding to the cyclic triads in~\eqref{eq:fano_triads}: all three outer sides and the inscribed circle $(1,2,4)$ are oriented counterclockwise, and all three medians pass inward through $e_7$. Any triple on a common line generates an associative quaternionic subalgebra ($\|\calA_3\| = 1$), whereas any three non-collinear points satisfy $\calA_3 = 0$. (b)~Root-split factorization of Theorem~\ref{thm:bilinear_dp}: bilinearity collapses the sum over $C_{N-1}$ trees into an exact $O(N^3)$ dynamic program weighted by the Catalan split probability $w_{L,m}$.}
\label{fig:fano_dp}
\end{figure}

\begin{proposition}[Sub-Unitary Contraction and Equality Criterion]
\label{prop:contraction}
In $\OO$ (or any normed division algebra), for all $x_1, \dots, x_N \in \OO$,
\begin{equation}
\label{eq:contraction}
\|\calA_N(x_1, \dots, x_N)\| \;\le\; \prod_{i=1}^N \|x_i\|.
\end{equation}
For nonzero $x_1, \dots, x_N \in \OO$, equality holds in~\eqref{eq:contraction} if and only if all $C_{N-1}$ tree products $\mathrm{prod}_T(x_1, \dots, x_N)$ ($T \in \calB_N$) are identical.
\end{proposition}

\begin{proof}
By~\eqref{eq:norm_mult}, every tree $T \in \calB_N$ satisfies $\|\mathrm{prod}_T(x_1, \dots, x_N)\| = \prod_{i=1}^N \|x_i\| =: R$. The triangle inequality gives $\|\calA_N\| \le C_{N-1}^{-1}\sum_{T \in \calB_N} \|\mathrm{prod}_T\| = R$. Because the Euclidean norm on $\OO \cong \R^8$ is strictly convex, an average of vectors of common radius $R > 0$ attains norm $R$ if and only if all vectors coincide.
\end{proof}

%==============================================================================
\section{The Universal Quaternionic Rigidity Theorem}
\label{sec:rigidity}
%==============================================================================

Proposition~\ref{prop:contraction} prompts a natural structural question: can a non-associative octonionic word ever preserve full norm $\|\calA_N\| = \prod \|x_i\|$ through hidden Moufang cancellations, repeated factors, or palindromic symmetries? We prove that the answer is strictly negative for all $N \ge 1$ and all nonzero octonions.

\begin{lemma}[Associator Kernel of a Quaternionic Subalgebra]
\label{lem:assoc_kernel}
Let $\HH_0 \subset \OO$ be a quaternionic subalgebra, and let $u, v \in \HH_0$ be two elements that do not lie in any common $2$-dimensional subalgebra $\R + \R a \subset \HH_0$ (equivalently, $\ImO(u)$ and $\ImO(v)$ are linearly independent in $\ImO(\HH_0)$). Then for any $z \in \OO$,
\begin{equation}
[u, v, z] = 0 \quad \Longleftrightarrow \quad z \in \HH_0.
\end{equation}
\end{lemma}

\begin{proof}
If $z \in \HH_0$, then $u, v, z \in \HH_0$ and associativity of $\HH_0$ gives $[u, v, z] = 0$. Conversely, choose a unit vector $\ell \in \HH_0^\perp$ so that $\OO = \HH_0 \oplus \HH_0\ell$ obeys the Cayley--Dickson rule~\eqref{eq:cayley_dickson}, and write $z = z_0 + q\ell$ with $z_0, q \in \HH_0$. Since the associator vanishes on real scalars and on $\HH_0$, we have $[u, v, z] = [u', v', q\ell]$, where $u' = \ImO(u)$ and $v' = \ImO(v)$ are linearly independent purely imaginary quaternions in $\ImO(\HH_0)$. Using~\eqref{eq:cayley_dickson} with $a = u', c = v' \in \ImO(\HH_0)$ and $b = d = 0$ on the first step, then multiplying by $0 + q\ell$:
\begin{align}
(u' v')(q\ell) &= \big(q(u' v')\big)\ell, \\
u'\big(v'(q\ell)\big) &= u'\big((qv')\ell\big) = \big((qv')u'\big)\ell = \big(q(v' u')\big)\ell,
\end{align}
where we used associativity of $\HH_0$ in the last step. Subtracting gives
\begin{equation}
[u, v, z] = \big(q[u', v']\big)\ell,
\end{equation}
where $[u', v'] = u' v' - v' u' \in \ImO(\HH_0)$ is the quaternionic commutator. Because $u'$ and $v'$ are linearly independent in $\ImO(\HH_0)$, their cross product $\frac{1}{2}[u', v']$ is nonzero. Since $\HH_0$ is a division algebra, $q[u', v'] = 0$ forces $q = 0$, hence $z = z_0 \in \HH_0$.
\end{proof}

\begin{theorem}[Universal Quaternionic Rigidity]
\label{thm:rigidity}
Let $N \ge 1$ and let $x_1, \dots, x_N \in \OO \setminus \{0\}$. Then
\begin{equation}
\|\calA_N(x_1, \dots, x_N)\| = \prod_{i=1}^N \|x_i\|
\end{equation}
if and only if the real subalgebra generated by $\{x_1, \dots, x_N\}$ is associative (equivalently, $\{x_1, \dots, x_N\}$ is contained in a single quaternionic subalgebra $\HH \subset \OO$).
\end{theorem}

\begin{proof}
The ``if'' direction is immediate since all $C_{N-1}$ parenthesizations coincide in any associative subalgebra. For the ``only if'' direction, assume $\|\calA_N(x_1, \dots, x_N)\| = \prod_{i=1}^N \|x_i\|$.

\medskip
\noindent\emph{Step 1 (Normalization and removal of real entries).}
Dividing each $x_i$ by $\|x_i\|$ preserves both the equality and the generated subalgebra, so assume $\|x_i\| = 1$. Any real entry $x_i \in \{\pm 1\}$ commutes and associates with all elements of $\OO$ and factors out identically from every $\mathrm{prod}_T(x_1, \dots, x_N)$. Deleting all real entries leaves $\|\calA_N\|$ unchanged and does not alter whether the remaining entries lie in a quaternionic subalgebra. Thus we may assume without loss of generality that $x_i \in \OO \setminus \R$ for every $i = 1, \dots, N$.

\medskip
\noindent\emph{Step 2 (Subword well-definedness and block associativity).}
By Proposition~\ref{prop:contraction}, all $C_{N-1}$ tree products $\mathrm{prod}_T(x_1, \dots, x_N)$ are equal. Consider any contiguous subword $w_{a,b} := (x_{a+1}, \dots, x_b)$ ($0 \le a < b \le N$). Fixing any outer parenthesization around the block $w_{a,b}$ defines a map $\Phi: \OO \to \OO$ formed by composing left and right multiplications by nonzero octonions. Because $\OO$ is a division algebra, $\Phi$ is injective. Therefore any two trees $T_1, T_2 \in \calB_{b-a}$ must satisfy $\mathrm{prod}_{T_1}(w_{a,b}) = \mathrm{prod}_{T_2}(w_{a,b})$, since otherwise embedding $T_1$ and $T_2$ inside $\Phi$ would produce two distinct values of $\mathrm{prod}_T(x_1, \dots, x_N)$. Hence every contiguous subword $w_{a,b}$ has a uniquely defined nonzero product $P(a, b) \in \OO \setminus \{0\}$. Moreover, for any $0 \le a < b < c < d \le N$, comparing the two trees grouping $((w_{a,b} w_{b,c}) w_{c,d})$ and $(w_{a,b} (w_{b,c} w_{c,d}))$ inside the same outer tree gives
\begin{equation}
\label{eq:block_assoc}
[P(a, b),\, P(b, c),\, P(c, d)] = 0.
\end{equation}

\medskip
\noindent\emph{Step 3 (Minimal escape index from a quaternionic subalgebra).}
Suppose for contradiction that $\{x_1, \dots, x_N\}$ is not contained in any quaternionic subalgebra of $\OO$. Let $k \in \{3, \dots, N\}$ be the smallest index such that $\{x_1, \dots, x_k\}$ is not contained in any quaternionic subalgebra (note $k \ge 3$ by Artin's theorem). Then $\{x_1, \dots, x_{k-1}\}$ is contained in some quaternionic subalgebra $\HH_0 \subset \OO$. Furthermore, $\{x_1, \dots, x_{k-1}\}$ cannot lie inside any $2$-dimensional complex subfield $\C' = \R + \R y \subset \HH_0$: if it did, then $\{x_1, \dots, x_k\}$ would lie in the subalgebra generated by the two elements $\{y, x_k\}$, which is associative (hence contained in a quaternionic subalgebra) by Artin's theorem, contradicting the minimality of $k$. Thus $\{x_1, \dots, x_{k-1}\}$ generates $\HH_0$, and $x_k \notin \HH_0$.

\medskip
\noindent\emph{Step 4 (Last escape index from a $2$-dimensional subfield of $\HH_0$).}
Since $x_{k-1} \in \HH_0 \setminus \R$, the singleton $\{x_{k-1}\}$ lies in the unique $2$-dimensional complex subfield $\C_1 := \R + \R x_{k-1} \subset \HH_0$. Because $\{x_1, \dots, x_{k-1}\} \not\subset \C_1$, there exists a largest index $m \in \{1, \dots, k-2\}$ such that $x_m \notin \C_1$. Thus
\begin{equation}
x_m \in \HH_0 \setminus \C_1, \qquad x_{m+1}, \dots, x_{k-1} \in \C_1 \setminus \R.
\end{equation}
Since $\C_1$ is closed under multiplication, $P(m, k-1) \in \C_1$. We consider two exhaustive cases:
\begin{itemize}
    \item \textbf{Case 1: $P(m, k-1) \notin \R$.} Set $u := P(m-1, m) = x_m \in \HH_0 \setminus \C_1$ and $v := P(m, k-1) \in \C_1 \setminus \R$. Then $u$ and $v$ do not lie in a common $2$-dimensional subfield of $\HH_0$ (any $2$-dimensional subfield containing $v \in \C_1 \setminus \R$ must equal $\C_1$, which excludes $u$). By~\eqref{eq:block_assoc} with $(a,b,c,d) = (m-1, m, k-1, k)$, we obtain $[u, v, x_k] = 0$.
    \item \textbf{Case 2: $P(m, k-1) \in \R$.} Since $x_{k-1} \notin \R$, we must have $m + 1 < k - 1$, so $P(m+1, k-1)$ is a well-defined element of $\C_1$. From $P(m, k-1) = x_{m+1} P(m+1, k-1) \in \R \setminus \{0\}$ and $x_{m+1} \in \C_1 \setminus \R$, it follows that $v := P(m+1, k-1) = x_{m+1}^{-1} P(m, k-1) \in \C_1 \setminus \R$. Meanwhile, $u := P(m-1, m+1) = x_m x_{m+1}$ lies in $\HH_0 \setminus \C_1$ (if $x_m x_{m+1} \in \C_1$, then multiplying on the right by $x_{m+1}^{-1} \in \C_1$ would force $x_m \in \C_1$, a contradiction). Applying~\eqref{eq:block_assoc} with $(a,b,c,d) = (m-1, m+1, k-1, k)$ again gives $[u, v, x_k] = 0$ with $u \in \HH_0 \setminus \C_1$ and $v \in \C_1 \setminus \R$.
\end{itemize}
In both cases, $u, v \in \HH_0$ have linearly independent imaginary parts in $\ImO(\HH_0)$ and satisfy $[u, v, x_k] = 0$. By Lemma~\ref{lem:assoc_kernel}, this forces $x_k \in \HH_0$, contradicting Step~3.
\end{proof}

\subsection{Fano Basis Words: $\GL(3, \F_2)$ Quantization and Exact Parity Obstruction}

Restricting to words $w = (e_{i_1}, \dots, e_{i_N}) \in \{e_1, \dots, e_7\}^N$ in the imaginary Fano basis reveals sharp finite-geometric structure.

\begin{proposition}[Rational Tree Sums and $\GL(3,\F_2)$ Orbit Quantization]
\label{prop:fano_orbit}
For any Fano basis word $w = (e_{i_1}, \dots, e_{i_N}) \in \{e_1, \dots, e_7\}^N$, there exist a unique basis index $\mu(w) \in \{0, 1, \dots, 7\}$ (independent of parenthesization) and an integer signed tree count $s(w) \in \{-C_{N-1}, \dots, C_{N-1}\}$ such that
\begin{equation}
\calA_N(w) = \frac{s(w)}{C_{N-1}}\,e_{\mu(w)}, \qquad s(w) \equiv C_{N-1} \pmod{2}.
\end{equation}
Moreover, the number of non-quaternionic words $w \in \{e_1, \dots, e_7\}^N$ achieving any fixed value $|s(w)| = m$ is an exact integer multiple of $168 = |\GL(3, \F_2)|$.
\end{proposition}

\begin{proof}
Modulo signs $\{\pm 1\}$, the octonionic basis loop $Q_\OO = \{\pm e_0, \dots, \pm e_7\}$ projects onto the elementary abelian group $Q_\OO / \{\pm 1\} \cong \F_2^3$ (identifying $e_0$ with $0 \in \F_2^3$ and the seven Fano lines with the seven $2$-dimensional subspaces of $\F_2^3$). Hence every parenthesization $\mathrm{prod}_T(w)$ equals $\sigma_T(w)\,e_{\mu(w)}$ with $\sigma_T(w) \in \{\pm 1\}$ and $\mu(w)$ independent of $T$. Thus $s(w) = \sum_{T \in \calB_N} \sigma_T(w) \equiv C_{N-1} \pmod{2}$.

For the orbit divisibility, every automorphism $\pi \in \Aut(\mathbb{PG}(2,\F_2)) \cong \GL(3,\F_2)$ (of order $168$) lifts to a signed basis permutation $\tilde{\pi} \in \Aut(\OO) = G_2$, where $\tilde{\pi}(e_a) = \eta_a e_{\pi(a)}$ with $\eta_a \in \{\pm 1\}$. Because $\tilde{\pi}$ is an algebra automorphism, $\tilde{\pi}(\calA_N(w)) = \calA_N(\tilde{\pi}(w))$, whence $|s(\pi(w))| = |s(w)|$. If $w$ is non-quaternionic, its letters contain three non-collinear Fano points (a basis of $\F_2^3$). Since any element of $\GL(3,\F_2)$ fixing a basis of $\F_2^3$ is the identity, $\GL(3,\F_2)$ acts freely on the set of non-quaternionic words, so every orbit has exact size $168$.
\end{proof}

\begin{corollary}[Dyadic Parity Obstruction to Exact Cancellation]
\label{cor:parity}
Exact cancellation $\calA_N(w) = 0$ is impossible for any Fano basis word $w \in \{e_1, \dots, e_7\}^N$ whenever $N = 2^k$ ($k \ge 1$); in fact, $|\calA_{2^k}(w)| \ge C_{2^k - 1}^{-1}$ for all $7^{2^k}$ words.
\end{corollary}

\begin{proof}
By Legendre's formula, $v_2(\binom{2m}{m}) = 2 s_2(m) - s_2(2m) = s_2(m)$, where $s_2(m)$ is the binary digit sum of $m$. Writing $m + 1 = 2^r q$ with $q$ odd ($r = v_2(m+1)$) gives $m = 2^r(q-1) + (2^r - 1)$, whose binary expansions do not overlap, so $s_2(m) = s_2(q-1) + r$ and hence $v_2(C_m) = s_2(m) - v_2(m+1) = s_2(q-1)$. Thus $C_m$ is odd if and only if $q = 1$, i.e., $m + 1 = 2^r \iff N = 2^r$. By Proposition~\ref{prop:fano_orbit}, $s(w) \equiv C_{N-1} \equiv 1 \pmod{2}$ whenever $N = 2^k$, so $s(w) \neq 0$ and $|s(w)| \ge 1$.
\end{proof}

\subsection{Exact Exhaustive Fano Distributions through $N = 6$}

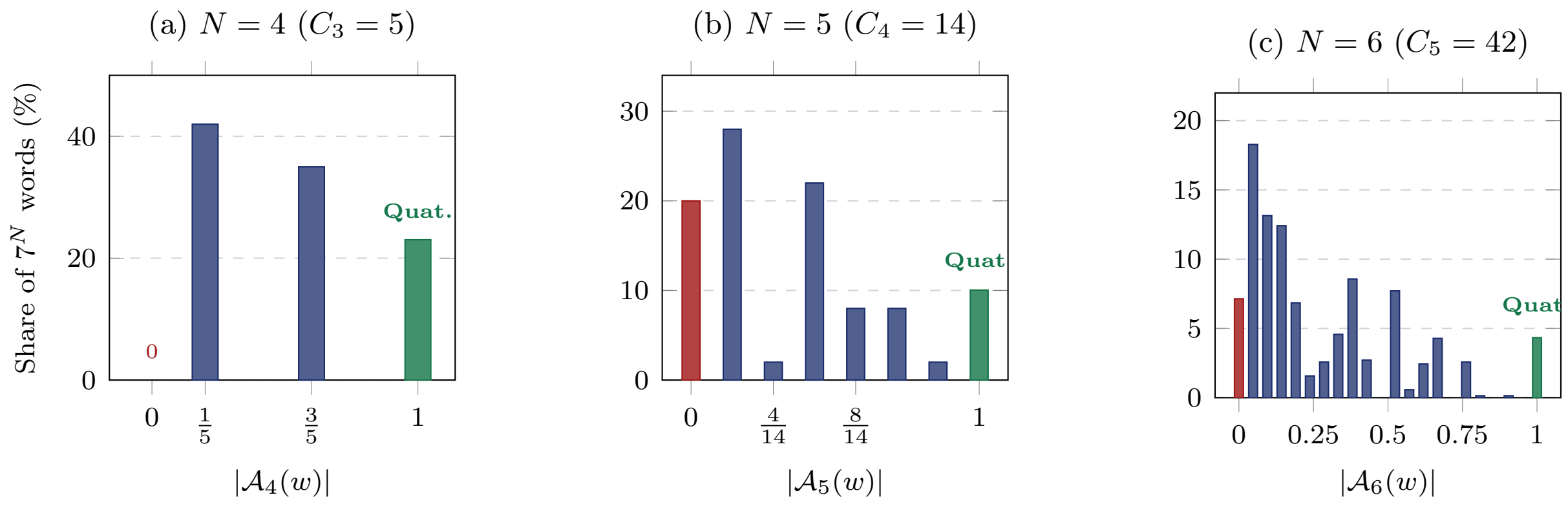
\begin{figure}[htbp]
\centering
\begin{tikzpicture}
\begin{axis}[
    width=0.32\textwidth,
    height=4.7cm,
    ybar,
    bar width=7.5pt,
    bar shift=0pt,
    xmin=-0.16, xmax=1.14,
    ymin=0, ymax=50,
    xtick={0,0.2,0.6,1.0},
    xticklabels={$0$,$\frac{1}{5}$,$\frac{3}{5}$,$1$},
    xlabel={$|\calA_4(w)|$},
    ylabel={Share of $7^N$ words (\%)},
    title={\footnotesize (a) $N=4$ ($C_3=5$)},
    tick label style={font=\scriptsize},
    label style={font=\scriptsize},
    ymajorgrids=true,
    grid style={dashed, slateink!25}
]
    \addplot[fill=navy!75, draw=navy] coordinates {(0.2, 41.98) (0.6, 34.99)};
    \addplot[fill=forest!80, draw=forest] coordinates {(1.0, 23.03)};
    \node[font=\tiny\bfseries, text=darkred, anchor=south] at (axis cs:0.0, 1.5) {$0$};
    \node[font=\tiny\bfseries, text=forest, anchor=south] at (axis cs:1.0, 23.8) {Quat.};
\end{axis}
\end{tikzpicture}
\hfill
\begin{tikzpicture}
\begin{axis}[
    width=0.32\textwidth,
    height=4.7cm,
    ybar,
    bar width=5.2pt,
    bar shift=0pt,
    xmin=-0.10, xmax=1.10,
    ymin=0, ymax=34,
    xtick={0,0.2857,0.5714,1.0},
    xticklabels={$0$,$\frac{4}{14}$,$\frac{8}{14}$,$1$},
    xlabel={$|\calA_5(w)|$},
    title={\footnotesize (b) $N=5$ ($C_4=14$)},
    tick label style={font=\scriptsize},
    label style={font=\scriptsize},
    ymajorgrids=true,
    grid style={dashed, slateink!25}
]
    \addplot[fill=darkred!80, draw=darkred] coordinates {(0.0, 19.99)};
    \addplot[fill=navy!75, draw=navy] coordinates {
        (0.1429, 27.99) (0.2857, 2.00) (0.4286, 21.99)
        (0.5714, 8.00) (0.7143, 8.00) (0.8571, 2.00)
    };
    \addplot[fill=forest!80, draw=forest] coordinates {(1.0, 10.04)};
    \node[font=\tiny\bfseries, text=forest, anchor=south] at (axis cs:1.0, 10.7) {Quat.};
\end{axis}
\end{tikzpicture}
\hfill
\begin{tikzpicture}
\begin{axis}[
    width=0.32\textwidth,
    height=4.7cm,
    ybar,
    bar width=2.6pt,
    bar shift=0pt,
    xmin=-0.08, xmax=1.08,
    ymin=0, ymax=22,
    xtick={0,0.25,0.5,0.75,1.0},
    xticklabels={$0$,$0.25$,$0.5$,$0.75$,$1$},
    xlabel={$|\calA_6(w)|$},
    title={\footnotesize (c) $N=6$ ($C_5=42$)},
    tick label style={font=\scriptsize},
    label style={font=\scriptsize},
    ymajorgrids=true,
    grid style={dashed, slateink!25}
]
    \addplot[fill=darkred!80, draw=darkred] coordinates {(0.0, 7.14)};
    \addplot[fill=navy!75, draw=navy] coordinates {
        (0.0476, 18.28) (0.0952, 13.14) (0.1429, 12.42) (0.1905, 6.85)
        (0.2381, 1.57) (0.2857, 2.57) (0.3333, 4.57) (0.3810, 8.57)
        (0.4286, 2.71) (0.5238, 7.71) (0.5714, 0.57) (0.6190, 2.43)
        (0.6667, 4.28) (0.7619, 2.57) (0.8095, 0.14) (0.9048, 0.14)
    };
    \addplot[fill=forest!80, draw=forest] coordinates {(1.0, 4.33)};
    \node[font=\tiny\bfseries, text=forest, anchor=south] at (axis cs:1.0, 5.0) {Quat.};
\end{axis}
\end{tikzpicture}
\caption{Exact exhaustive distributions of $|\calA_N(w)| = |s(w)|/C_{N-1}$ over all $7^N$ Fano basis words at $N=4, 5, 6$. Crimson bars mark exact zero cancellations ($s=0$), blue bars mark non-quaternionic intermediate values (every count an exact multiple of $168 = |\GL(3,\F_2)|$, Proposition~\ref{prop:fano_orbit}), and green bars at $|\calA_N|=1$ consist exclusively of quaternionic words (Theorem~\ref{thm:rigidity}). Note the complete absence of zeros at $N=4=2^2$ (Corollary~\ref{cor:parity}) and the sharp spectral gap below $|\calA_6|=1$.}
\label{fig:fano_dist}
\end{figure}

Using Theorem~\ref{thm:bilinear_dp} in exact integer arithmetic (exploiting the $7$-fold transitivity on the first letter at $N = 6$), we enumerated all $7^N$ Fano basis words for $N = 3, 4, 5, 6$. Table~\ref{tab:fano_exact} and Figure~\ref{fig:fano_dist} display the exact distributions and moments. In all four cases, the count of $|s(w)| = C_{N-1}$ words equals the exact number of quaternionic words $7 + 21(2^N - 2) + 7(3^N - 3\cdot 2^N + 3)$ with zero exceptions, confirming Theorem~\ref{thm:rigidity}; every non-quaternionic count is divisible by $168$, and exact zero cancellations vanish at $N = 4 = 2^2$ while accounting for $48.98\%$ ($N=3$), $19.99\%$ ($N=5$), and $7.14\%$ ($N=6$) of all words. Notice also the clean spectral gap below $|\calA_6| = 1$ in Figure~\ref{fig:fano_dist}(c): the values $|s| \in \{20, 30, 36, 40\}$ never occur, and only a single $\GL(3,\F_2)$ orbit ($168$ words, $0.14\%$) reaches $|s| = 38$ ($|\calA_6| = 38/42$).

\begin{table}[htbp]
\centering
\caption{Exact exhaustive statistics of $\calA_N(w)$ over all $7^N$ Fano basis words for $N = 3, 4, 5, 6$.}
\label{tab:fano_exact}
\vspace{0.4em}
\small
\begin{tabular}{@{}rrrrrrl@{}}
\toprule
$N$ & $C_{N-1}$ & $7^N$ & Quat.\ ($|s|=C_{N-1}$) & Zeros ($s=0$) & Exact $\E[|\calA_N|^2]$ & Decimal \\
\midrule
$3$ & $2$  & $343$     & $175$   & $168 = 1 \times 168$   & $25/49$             & $0.51020408$ \\
$4$ & $5$  & $2{,}401$ & $553$   & $0$ (Corollary~\ref{cor:parity}) & $457/1225$ & $0.37306122$ \\
$5$ & $14$ & $16{,}807$ & $1{,}687$ & $3{,}360 = 20 \times 168$ & $27025/117649$ & $0.22970871$ \\
$6$ & $42$ & $117{,}649$ & $5{,}089$ & $8{,}400 = 50 \times 168$ & $16913/117649$ & $0.14375813$ \\
\bottomrule
\end{tabular}
\end{table}

%==============================================================================
\section{Exact $N = 3$ Associator Moments and Rigorous Exponential Decay}
\label{sec:exp_decay}
%==============================================================================

We now analyze $\calA_N(X_1, \dots, X_N)$ when $X_1, \dots, X_N$ are independent, identically distributed (i.i.d.) random unit octonions drawn from a probability measure $\mu$ on $S^7 \subset \OO$. Our four primary ensembles are:
\begin{enumerate}
    \item $\mu_{\mathrm{Fano}}$: uniform on the $7$ imaginary Fano basis units $\{e_1, \dots, e_7\}$;
    \item $\mu_{S^6}$: uniform (Haar/surface) measure on the unit imaginary sphere $S^6 \subset \ImO(\OO)$;
    \item $\mu_{V_\perp}$: uniform measure on the unit $5$-sphere in the $6$-dimensional transverse subspace $V_\perp := e_7^\perp = \Span_\R\{e_1, \dots, e_6\}$;
    \item $\mu_{S^7}$: uniform measure on the full unit octonion sphere $S^7 \subset \OO$.
\end{enumerate}

\subsection{Exact $N = 3$ Isotropic Associator Identity}

\begin{theorem}[Exact $N = 3$ Second Moment]
\label{thm:n3_exact}
Let $V = \Span_\R\{e_a : a \in I\} \subset \OO$ be a coordinate subspace of dimension $d = |I|$, and let $X_1, X_2, X_3 \in V \cap S^7$ be i.i.d.\ unit random octonions with isotropic second-moment covariance $\E[(X_1)_a (X_1)_b] = \frac{1}{d}\delta_{ab}$ for $a, b \in I$. Then
\begin{equation}
\label{eq:n3_general}
\E\big[\|\calA_3(X_1, X_2, X_3)\|^2\big] = 1 - \frac{1}{4}\E\big[\|[X_1, X_2, X_3]\|^2\big] = 1 - \frac{1}{4 d^3} \sum_{a,b,c \in I} \|[e_a, e_b, e_c]\|^2.
\end{equation}
In particular,
\begin{equation}
\label{eq:n3_values}
\E[\|\calA_3\|^2] = \begin{cases}
\dfrac{25}{49} \approx 0.510204 & \text{for } \mu_{\mathrm{Fano}} \text{ and } \mu_{S^6} \; (d = 7,\; \sum \|[e_a,e_b,e_c]\|^2 = 672), \\[0.6em]
\dfrac{5}{9} \approx 0.555556 & \text{for } \mu_{V_\perp} \; (d = 6,\; \sum \|[e_a,e_b,e_c]\|^2 = 384), \\[0.6em]
\dfrac{43}{64} = 0.671875 & \text{for } \mu_{S^7} \; (d = 8,\; \sum \|[e_a,e_b,e_c]\|^2 = 672).
\end{cases}
\end{equation}
\end{theorem}

\begin{proof}
For any three unit octonions $x, y, z \in S^7$, write $p := (xy)z$ and $q := x(yz)$, so $\|p\| = \|q\| = 1$. By the parallelogram identity,
\begin{equation}
\|\calA_3(x, y, z)\|^2 = \left\|\frac{p + q}{2}\right\|^2 = \frac{1}{2}(\|p\|^2 + \|q\|^2) - \frac{1}{4}\|p - q\|^2 = 1 - \frac{1}{4}\|[x, y, z]\|^2.
\end{equation}
Expanding $X_1 = \sum_{a \in I} (X_1)_a e_a$, $X_2 = \sum_{b \in I} (X_2)_b e_b$, $X_3 = \sum_{c \in I} (X_3)_c e_c$ and using independence together with $\E[(X_r)_a (X_r)_{a'}] = \frac{1}{d}\delta_{aa'}$ yields
\begin{align}
\E\big[\|[X_1, X_2, X_3]\|^2\big]
&= \sum_{a,b,c,a',b',c' \in I} \E[(X_1)_a(X_1)_{a'}]\,\E[(X_2)_b(X_2)_{b'}]\,\E[(X_3)_c(X_3)_{c'}]\,\langle [e_a,e_b,e_c], [e_{a'},e_{b'},e_{c'}] \rangle \nonumber\\
&= \frac{1}{d^3} \sum_{a,b,c \in I} \|[e_a, e_b, e_c]\|^2,
\end{align}
proving~\eqref{eq:n3_general}. For each $a \in \{1, \dots, 7\}$, there are $6$ choices of $b \neq a$ and $4$ choices of $c$ outside the unique Fano line through $a$ and $b$, each contributing $\|[e_a, e_b, e_c]\|^2 = 4$; hence $\sum_{b,c=1}^7 \|[e_a, e_b, e_c]\|^2 = 6 \times 4 \times 4 = 96$ (the octonionic associator Casimir~\cite{smith2026octonion}) and $\sum_{a,b,c=1}^7 \|[e_a, e_b, e_c]\|^2 = 7 \times 96 = 672$. Since $e_0 = 1$ has zero associator, the sum over $I = \{0, \dots, 7\}$ is also $672$. On $V_\perp = \Span_\R\{e_1, \dots, e_6\}$, for each $a \in \{1, \dots, 6\}$ the Fano line $(7, a, a')$ removes one partner $a' \in \{1,\dots,6\}$ (leaving $4$ choices of $b$ outside that line, each yielding $3$ choices of $c \in \{1,\dots,6\}$ off the line through $a, b$, plus $1$ choice $b = a'$ yielding $4$ choices of $c \in \{1,\dots,6\} \setminus \{a,a'\}$), giving $(4 \times 3 + 1 \times 4) \times 4 = 64$ per index and $6 \times 64 = 384$ in total. Substituting into~\eqref{eq:n3_general} gives $1 - \frac{672}{4 \cdot 343} = \frac{25}{49}$, $1 - \frac{384}{4 \cdot 216} = \frac{5}{9}$, and $1 - \frac{672}{4 \cdot 512} = \frac{43}{64}$.
\end{proof}

\subsection{Rigorous Proof of Exponential Decay via Generating-Function Singularity Analysis}

We now prove that for \emph{any} i.i.d.\ unit-octonion distribution $\mu$ not supported inside a single quaternionic subalgebra, $\E[\|\calA_N\|^2]$ decays at least exponentially in $N$ and $\|\calA_N\| \to 0$ exponentially almost surely.

\begin{lemma}[Catalan Majorant Inequality for $L^2$ Norms]
\label{lem:majorant}
Let $X_1, X_2, \dots \in S^7 \subset \OO$ be i.i.d.\ random unit octonions, and define
\begin{equation}
a_N := \Big(\E\big[\|\calA_N(X_1, \dots, X_N)\|^2\big]\Big)^{1/2} \in [0, 1].
\end{equation}
Then $a_1 = a_2 = 1$, and for every $N \ge 3$,
\begin{equation}
\label{eq:a_ineq}
a_N \;\le\; \sum_{k=1}^{N-1} w_{N, k}\,a_k\,a_{N-k}.
\end{equation}
\end{lemma}

\begin{proof}
By Theorem~\ref{thm:bilinear_dp}, $\calA_N = \sum_{k=1}^{N-1} w_{N,k}\,\calA(0,k)\,\calA(k,N)$. Applying Minkowski's triangle inequality in $L^2(\Prob; \OO)$ gives
\begin{equation}
a_N = \|\calA_N\|_{L^2} \;\le\; \sum_{k=1}^{N-1} w_{N,k}\,\big\|\calA(0,k)\,\calA(k,N)\big\|_{L^2}.
\end{equation}
For each fixed $k \in \{1, \dots, N-1\}$, $\calA(0, k)$ is a function of $(X_1, \dots, X_k)$ and $\calA(k, N)$ is a function of $(X_{k+1}, \dots, X_N)$, so these two random octonions are independent. By exact norm multiplicativity~\eqref{eq:norm_mult} and stationarity,
\begin{equation}
\E\big[\|\calA(0,k)\,\calA(k,N)\|^2\big] = \E\big[\|\calA(0,k)\|^2\,\|\calA(k,N)\|^2\big] = \E\big[\|\calA_k\|^2\big]\,\E\big[\|\calA_{N-k}\|^2\big] = a_k^2\,a_{N-k}^2,
\end{equation}
taking square roots yields~\eqref{eq:a_ineq}.
\end{proof}

\begin{theorem}[Rigorous Exponential Decay and Almost-Sure Convergence]
\label{thm:exp_decay}
Let $X_1, X_2, \dots \in S^7 \subset \OO$ be i.i.d.\ random unit octonions with $a_3 = (\E[\|\calA_3\|^2])^{1/2} < 1$, and set $\delta := 2(1 - a_3) \in (0, 2]$. Let $u_0 \in (1, \frac{3}{2})$ be the unique root in $(1, \frac{3}{2})$ of the cubic equation
\begin{equation}
\label{eq:cubic_u0}
1 - u + \frac{\delta}{16}\,u^3 = 0.
\end{equation}
Then:
\begin{enumerate}
    \item[(i)] \textbf{Explicit $L^2$ Exponential Decay:} The sequence $a_N$ is bounded above by the exact Catalan-defect solution $a_N^*$, which satisfies the sharp asymptotic
    \begin{equation}
    \label{eq:sharp_upper}
    a_N \;\le\; a_N^* \;=\; \sqrt{3 - 2 u_0}\;u_0^{-N}\,\big(1 + O(N^{-1})\big) \qquad (N \to \infty).
    \end{equation}
    In particular, there exists an explicit constant $M > 0$ such that $\E[\|\calA_N\|^2] \le M\,u_0^{-2N}$ for all $N \ge 1$.
    \item[(ii)] \textbf{Almost-Sure Exponential Convergence:} Almost surely, $\|\calA_N(X_1, \dots, X_N)\| \to 0$ exponentially fast, with
    \begin{equation}
    \label{eq:as_lyap_bound}
    \liminf_{N \to \infty} \left(-\frac{1}{N}\log \|\calA_N(X_1, \dots, X_N)\|\right) \;\ge\; \log u_0 \;>\; 0 \qquad \text{almost surely.}
    \end{equation}
    \item[(iii)] \textbf{Multi-Order Tightening ($K \ge 3$):} Given exact values (or upper bounds) $a_1, \dots, a_K$ with nonnegative defects $d_m := \sum_{k=1}^{m-1} C_{k-1}C_{m-k-1}a_k a_{m-k} - C_{m-1}a_m \ge 0$ ($3 \le m \le K$), let $u_0^{(K)} > 1$ be the smallest positive root of
    \begin{equation}
    \label{eq:poly_PK}
    1 - u + 4\sum_{m=3}^K d_m \left(\frac{u}{4}\right)^m = 0.
    \end{equation}
    Then $\limsup_{N \to \infty} a_N^{1/N} \le (u_0^{(K)})^{-1}$, and~\eqref{eq:as_lyap_bound} holds with $\log u_0^{(K)}$ in place of $\log u_0$.
\end{enumerate}
\end{theorem}

\begin{proof}
Define $b_N := C_{N-1} a_N$ for $N \ge 1$. Then $b_1 = 1$, $b_2 = 1$, $b_3 = 2 a_3 = 2 - \delta < 2$, and Lemma~\ref{lem:majorant} translates to
\begin{equation}
b_N \;\le\; \sum_{k=1}^{N-1} b_k\,b_{N-k} \qquad (N \ge 3).
\end{equation}
For any fixed order $K \ge 3$, define the majorant sequence $\{b_N^*\}_{N \ge 1}$ by $b_m^* := b_m$ for $1 \le m \le K$ and $b_N^* := \sum_{k=1}^{N-1} b_k^* b_{N-k}^*$ for $N > K$. Since $b_m \ge 0$, strong induction gives $0 \le b_N \le b_N^* \le C_{N-1}$ for all $N \ge 1$.

Form the ordinary generating function $B^*(z) := \sum_{N=1}^\infty b_N^* z^N$, which converges absolutely on $|z| < 1/4$ because $0 \le b_N^* \le C_{N-1}$. Comparing coefficients of $z^N$ on both sides shows that $B^*(z)$ satisfies the quadratic functional equation
\begin{equation}
B^*(z)^2 - B^*(z) + z - \sum_{m=3}^K d_m z^m = 0,
\end{equation}
where $d_m = \sum_{k=1}^{m-1} b_k^* b_{m-k}^* - b_m^* \ge 0$ (guaranteed by Lemma~\ref{lem:majorant}), with $d_3 = \delta = 2(1 - a_3) > 0$. Solving for the branch vanishing at $z = 0$ gives
\begin{equation}
B^*(z) = \frac{1 - \sqrt{P_K(z)}}{2}, \qquad P_K(z) := 1 - 4z + 4\sum_{m=3}^K d_m z^m.
\end{equation}

\medskip
\noindent\emph{Proof of (i) ($K = 3$).}
Setting $u = 4z$ and $c := \delta / 16 \in (0, \frac{1}{8}]$, we have $P_3(u/4) = f(u) := 1 - u + c u^3$. Observe that $f(0) = 1 > 0$, $f(1) = c > 0$, and
\begin{equation}
f\!\left(\frac{3}{2}\right) = -\frac{1}{2} + \frac{27}{8} c \;\le\; -\frac{1}{2} + \frac{27}{64} = -\frac{5}{64} < 0,
\end{equation}
while $\lim_{u \to \pm\infty} f(u) = \pm\infty$. Hence $f$ has three distinct real roots satisfying $u_2 < 0 < 1 < u_0 < \frac{3}{2} < u_1$. Moreover, because $c u_0^3 = u_0 - 1$, evaluating at $-u_0$ gives
\begin{equation}
f(-u_0) = 1 + u_0 - c u_0^3 = 1 + u_0 - (u_0 - 1) = 2 > 0,
\end{equation}
which forces $u_2 < -u_0$, i.e., $|u_2| > u_0$. Therefore $z_0 := u_0 / 4$ is the \emph{unique} singularity of $B^*(z)$ of smallest modulus in $\C$, and it is a simple square-root branch point with
\begin{equation}
-z_0 P_3'(z_0) = u_0 - 3 c u_0^3 = u_0 - 3(u_0 - 1) = 3 - 2 u_0 > 0.
\end{equation}
By the Flajolet--Sedgewick Transfer Theorem for algebraic singularities~\cite[Theorem~VI.1]{flajolet2009},
\begin{equation}
b_N^* = [z^N]\,\frac{1 - \sqrt{P_3(z)}}{2} = \frac{\sqrt{-z_0 P_3'(z_0)}}{4\sqrt{\pi}}\,N^{-3/2}\,z_0^{-N}\,\big(1 + O(N^{-1})\big).
\end{equation}
Dividing by $C_{N-1} = \frac{1}{4\sqrt{\pi}} N^{-3/2} 4^N (1 + O(N^{-1}))$ cancels the $N^{-3/2}/(4\sqrt{\pi})$ prefactor identically and yields~\eqref{eq:sharp_upper}:
\begin{equation}
a_N^* = \frac{b_N^*}{C_{N-1}} = \sqrt{3 - 2 u_0}\;u_0^{-N}\,\big(1 + O(N^{-1})\big).
\end{equation}

\medskip
\noindent\emph{Proof of (iii) (general $K \ge 3$).}
Because $b_N^* = \sum_{k=1}^{N-1} b_k^* b_{N-k}^* \ge 2 b_1^* b_{N-1}^* = 2 b_{N-1}^* > 0$ for all $N > K$, the power series $B^*(z)$ has strictly positive coefficients and finite radius of convergence $R_K \le 1/2$. By Pringsheim's Theorem~\cite[Theorem~IV.6]{flajolet2009}, $z = R_K$ is a singularity of $B^*(z)$ on the positive real axis. Since $B^*(z) = \frac{1}{2}(1 - \sqrt{P_K(z)})$ is analytic away from roots of $P_K$ and $P_K''(z) = 4\sum_{m=3}^K m(m-1)d_m z^{m-2} > 0$ makes $P_K$ strictly convex on $(0,\infty)$ (so any double root would leave $\sqrt{P_K(z)}$ analytic there), $R_K = z_0^{(K)} = u_0^{(K)}/4$ must be the smallest positive root of $P_K(z) = 0$ and is necessarily simple. Since $P_K(1/4) = 4\sum_{m=3}^K d_m 4^{-m} > 0$, we have $u_0^{(K)} > 1$. By the Cauchy--Hadamard formula, $\limsup_{N\to\infty} (b_N^*)^{1/N} = 1/z_0^{(K)}$, and dividing by $C_{N-1}^{1/N} \to 4$ gives $\limsup_{N\to\infty} a_N^{1/N} \le (u_0^{(K)})^{-1}$.

\medskip
\noindent\emph{Proof of (ii) (almost-sure exponential convergence).}
Fix any $\varepsilon \in (0, \log u_0^{(K)})$ and set $\gamma := \log u_0^{(K)} - \varepsilon > 0$. By Markov's inequality and part~(iii) (or part~(i) for $K=3$), for sufficiently large $N$,
\begin{equation}
\Prob\big(\|\calA_N\| > e^{-\gamma N}\big) \;\le\; e^{2\gamma N}\,\E[\|\calA_N\|^2] = e^{2\gamma N}\,a_N^2 \;\le\; C_\varepsilon\,e^{-\varepsilon N}.
\end{equation}
Since $\sum_{N=1}^\infty e^{-\varepsilon N} < \infty$, the First Borel--Cantelli Lemma implies that almost surely $\|\calA_N\| \le e^{-\gamma N}$ for all sufficiently large $N$, hence $\liminf_{N\to\infty} (-\frac{1}{N}\log\|\calA_N\|) \ge \log u_0^{(K)} - \varepsilon$ almost surely. Taking a countable sequence $\varepsilon_j \downarrow 0$ yields~\eqref{eq:as_lyap_bound}.
\end{proof}

\begin{remark}[Strict Positivity for Arbitrary Non-Quaternionic Measures]
By Theorem~\ref{thm:rigidity}, if $\mu$ is \emph{any} Borel probability measure on $S^7$ whose support is not contained in a single quaternionic subalgebra $\HH \subset \OO$, then $\Prob(\|\calA_3(X_1, X_2, X_3)\| < 1) > 0$, which automatically implies $a_3 < 1$ and triggers Theorem~\ref{thm:exp_decay}. Thus exponential decay of $\calA_N$ is a universal phenomenon of non-quaternionic probability measures on $S^7$.
\end{remark}

Table~\ref{tab:rigorous_bounds} lists the explicit analytic bounds yielded by Theorems~\ref{thm:n3_exact} and~\ref{thm:exp_decay}, together with direct recurrence verification at $N = 150$.

\begin{table}[htbp]
\centering
\caption{Rigorous analytic upper bounds on $a_N = (\E[\|\calA_N\|^2])^{1/2}$ and lower bounds $\lambda_{\min} = \log u_0^{(K)}$ on the almost-sure exponential decay rate (Theorem~\ref{thm:exp_decay}).}
\label{tab:rigorous_bounds}
\vspace{0.4em}
\small
\begin{tabular}{@{}lcccccc@{}}
\toprule
Ensemble / Order $K$ & $a_3$ & $u_0^{(K)}$ & $1/u_0^{(K)}$ & $\lambda_{\min} = \log u_0^{(K)}$ & $\sqrt{3 - 2u_0}$ & $a_{150}^* u_0^{150}$ \\
\midrule
$S^6$ and Fano ($K = 3$) & $5/7$            & $1.040196$ & $0.961357$ & $0.039410$ & $0.958961$ & $0.958322$ \\
$V_\perp = e_7^\perp$ ($K = 3$) & $\sqrt{5}/3$ & $1.035324$ & $0.965881$ & $0.034715$ & $0.964029$ & $0.963470$ \\
$S^7$ ($K = 3$)          & $\sqrt{43}/8$    & $1.024218$ & $0.976355$ & $0.023929$ & $0.975482$ & $0.975103$ \\
\midrule
Fano tightened ($K = 4$) & $5/7$ & $1.058026$ & $0.945157$ & $0.056404$ & --- & --- \\
Fano tightened ($K = 5$) & $5/7$ & $1.073378$ & $0.931638$ & $0.070811$ & --- & --- \\
Fano tightened ($K = 6$) & $5/7$ & $1.085355$ & $0.921358$ & $0.081907$ & --- & --- \\
\bottomrule
\end{tabular}
\end{table}

%==============================================================================
\section{Numerical Measurement of Non-Associative Lyapunov Exponents}
\label{sec:numerics}
%==============================================================================

While Theorem~\ref{thm:exp_decay} proves strictly positive lower bounds $\lambda_{\min} > 0$, Minkowski's triangle inequality ignores destructive angular cancellations between different split indices $k$ in~\eqref{eq:dp_recurrence}, so the true decay rates are substantially faster than $\lambda_{\min}$. Using the exact $O(N^3)$ normalized recurrence of Theorem~\ref{thm:bilinear_dp}, we measured the empirical finite-$N$ quenched Lyapunov exponent
\begin{equation}
\label{eq:emp_lyap}
\hat{\lambda}_N := -\frac{1}{N}\,\E\big[\log \|\calA_N(X_1, \dots, X_N)\|\big]
\end{equation}
(conditioned on $\calA_N \neq 0$ for discrete Fano words, where exact zero cancellations occur with vanishingly small frequency for large $N$; only $1$ zero occurred in $400$ samples at $N = 10$ and $0$ in all samples for $N \ge 20$). Table~\ref{tab:lyapunov_num} reports the results across $N = 10$ to $N = 80$ with standard errors of the mean.

\begin{table}[htbp]
\centering
\caption{Exact-$O(N^3)$ Monte Carlo measurements of the quenched non-associative Lyapunov exponent $\hat{\lambda}_N = -\frac{1}{N}\E[\log \|\calA_N\|]$ ($\pm$ standard error) over $M$ independent realizations up to $N = 80$ ($C_{79} \approx 2.89 \times 10^{44}$ trees per realization).}
\label{tab:lyapunov_num}
\vspace{0.4em}
\small
\begin{tabular}{@{}rrcccc@{}}
\toprule
$N$ & $M$ & Fano $\hat{\lambda}_N$ & $S^6$ $\hat{\lambda}_N$ & $V_\perp$ $\hat{\lambda}_N$ & $S^7$ $\hat{\lambda}_N$ \\
\midrule
$10$ & $400$ & $0.3417 \pm 0.0115$ & $0.2100 \pm 0.0023$ & $0.1912 \pm 0.0022$ & $0.1696 \pm 0.0020$ \\
$20$ & $200$ & $0.3466 \pm 0.0063$ & $0.2429 \pm 0.0021$ & $0.2205 \pm 0.0025$ & $0.1925 \pm 0.0021$ \\
$30$ & $120$ & $0.3669 \pm 0.0067$ & $0.2537 \pm 0.0026$ & $0.2244 \pm 0.0026$ & $0.2083 \pm 0.0023$ \\
$40$ &  $80$ & $0.3639 \pm 0.0078$ & $0.2628 \pm 0.0024$ & $0.2340 \pm 0.0023$ & $0.2151 \pm 0.0020$ \\
$60$ &  $50$ & $0.3490 \pm 0.0076$ & $0.2650 \pm 0.0026$ & $0.2397 \pm 0.0025$ & $0.2170 \pm 0.0018$ \\
$80$ &  $30$ & $0.3439 \pm 0.0063$ & $0.2679 \pm 0.0027$ & $0.2443 \pm 0.0028$ & $0.2217 \pm 0.0028$ \\
\midrule
\multicolumn{2}{@{}l}{Effective $e^{-\hat{\lambda}_{80}}$} & $0.709$ & $0.765$ & $0.783$ & $0.801$ \\
\bottomrule
\end{tabular}
\end{table}

\begin{figure}[htbp]
\centering
\begin{tikzpicture}
\begin{axis}[
    width=0.48\textwidth,
    height=5.6cm,
    xmin=0, xmax=84,
    ymin=-1.2, ymax=31,
    xlabel={Word length $N$},
    ylabel={$\E[-\log \|\calA_N\|]$},
    title={\footnotesize (a) Linear growth of $-\E[\log\|\calA_N\|]$ through $N=80$},
    tick label style={font=\scriptsize},
    label style={font=\scriptsize},
    legend pos=north west,
    legend style={font=\scriptsize, fill=white, fill opacity=0.88, draw=slateink!40, cells={anchor=west}},
    grid=major,
    grid style={dashed, slateink!25}
]
    \addplot[thick, darkred, mark=*, mark size=1.8pt, error bars/.cd, y dir=both, y explicit]
        coordinates {
            (0,0) +- (0,0) (10,3.417) +- (0,0.115) (20,6.932) +- (0,0.126)
            (30,11.007) +- (0,0.201) (40,14.556) +- (0,0.312)
            (60,20.940) +- (0,0.456) (80,27.512) +- (0,0.504)
        };
    \addlegendentry{Fano ($\hat{\lambda}_{80}=0.344$)}

    \addplot[thick, navy, mark=square*, mark size=1.7pt, error bars/.cd, y dir=both, y explicit]
        coordinates {
            (0,0) +- (0,0) (10,2.100) +- (0,0.023) (20,4.858) +- (0,0.042)
            (30,7.611) +- (0,0.078) (40,10.512) +- (0,0.096)
            (60,15.900) +- (0,0.156) (80,21.432) +- (0,0.216)
        };
    \addlegendentry{$S^6$ ($\hat{\lambda}_{80}=0.268$)}

    \addplot[thick, tealdark, mark=triangle*, mark size=2.1pt, error bars/.cd, y dir=both, y explicit]
        coordinates {
            (0,0) +- (0,0) (10,1.912) +- (0,0.022) (20,4.410) +- (0,0.050)
            (30,6.732) +- (0,0.078) (40,9.360) +- (0,0.092)
            (60,14.382) +- (0,0.150) (80,19.544) +- (0,0.224)
        };
    \addlegendentry{$V_\perp$ ($\hat{\lambda}_{80}=0.244$)}

    \addplot[thick, amberdark, mark=diamond*, mark size=2.1pt, error bars/.cd, y dir=both, y explicit]
        coordinates {
            (0,0) +- (0,0) (10,1.696) +- (0,0.020) (20,3.850) +- (0,0.042)
            (30,6.249) +- (0,0.069) (40,8.604) +- (0,0.080)
            (60,13.020) +- (0,0.108) (80,17.736) +- (0,0.224)
        };
    \addlegendentry{$S^7$ ($\hat{\lambda}_{80}=0.222$)}

    \addplot[thick, forest, dashed] coordinates {(0,0) (80,0)};
    \addlegendentry{$\HH$ ($\lambda=0$, associative)}
\end{axis}
\end{tikzpicture}
\hfill
\begin{tikzpicture}
\begin{axis}[
    width=0.48\textwidth,
    height=5.6cm,
    xmin=5, xmax=85,
    ymin=0, ymax=0.43,
    xlabel={Word length $N$},
    ylabel={Decay exponent $\hat{\lambda}_N$ and bounds $\log u_0^{(K)}$},
    title={\footnotesize (b) Empirical $\hat{\lambda}_N$ vs.\ rigorous lower bounds},
    tick label style={font=\scriptsize},
    label style={font=\scriptsize},
    grid=major,
    grid style={dashed, slateink!25}
]
    \addplot[thick, darkred, mark=*, mark size=1.8pt, error bars/.cd, y dir=both, y explicit]
        coordinates {
            (10,0.3417) +- (0,0.0115) (20,0.3466) +- (0,0.0063)
            (30,0.3669) +- (0,0.0067) (40,0.3639) +- (0,0.0078)
            (60,0.3490) +- (0,0.0076) (80,0.3439) +- (0,0.0063)
        };
    \addplot[thick, navy, mark=square*, mark size=1.7pt, error bars/.cd, y dir=both, y explicit]
        coordinates {
            (10,0.2100) +- (0,0.0023) (20,0.2429) +- (0,0.0021)
            (30,0.2537) +- (0,0.0026) (40,0.2628) +- (0,0.0024)
            (60,0.2650) +- (0,0.0026) (80,0.2679) +- (0,0.0027)
        };
    \addplot[thick, tealdark, mark=triangle*, mark size=2.1pt, error bars/.cd, y dir=both, y explicit]
        coordinates {
            (10,0.1912) +- (0,0.0022) (20,0.2205) +- (0,0.0025)
            (30,0.2244) +- (0,0.0026) (40,0.2340) +- (0,0.0023)
            (60,0.2397) +- (0,0.0025) (80,0.2443) +- (0,0.0028)
        };
    \addplot[thick, amberdark, mark=diamond*, mark size=2.1pt, error bars/.cd, y dir=both, y explicit]
        coordinates {
            (10,0.1696) +- (0,0.0020) (20,0.1925) +- (0,0.0021)
            (30,0.2083) +- (0,0.0023) (40,0.2151) +- (0,0.0020)
            (60,0.2170) +- (0,0.0018) (80,0.2217) +- (0,0.0028)
        };

    % Rigorous analytic lower bounds from Theorem 5.2 (Table 2)
    \addplot[thick, darkred, densely dashed] coordinates {(5,0.0819) (85,0.0819)};
    \node[font=\tiny\bfseries, text=darkred, anchor=south west] at (axis cs:8, 0.084) {Fano bound $\log u_0^{(6)} = 0.0819$};

    \addplot[thick, navy, densely dotted] coordinates {(5,0.0394) (85,0.0394)};
    \node[font=\tiny, text=navy, anchor=south west] at (axis cs:8, 0.041) {$S^6$/Fano ($K=3$): $0.0394$};

    \addplot[thick, amberdark, dotted] coordinates {(5,0.0239) (85,0.0239)};
    \node[font=\tiny, text=amberdark, anchor=north east] at (axis cs:83, 0.023) {$S^7$ ($K=3$): $0.0239$};

    % Direct curve labels on right
    \node[font=\tiny\bfseries, text=darkred, anchor=south east] at (axis cs:82, 0.355) {Fano};
    \node[font=\tiny\bfseries, text=navy, anchor=south east] at (axis cs:82, 0.273) {$S^6$};
    \node[font=\tiny\bfseries, text=tealdark, anchor=west] at (axis cs:62, 0.246) {$V_\perp$};
    \node[font=\tiny\bfseries, text=amberdark, anchor=north east] at (axis cs:82, 0.214) {$S^7$};
\end{axis}
\end{tikzpicture}
\caption{Exact-$O(N^3)$ non-associative decay through $N = 80$ ($C_{79} \approx 2.89 \times 10^{44}$ trees per realization). (a)~Mean negative log-norm $\E[-\log\|\calA_N\|]$ grows linearly in $N$ across all four non-associative ensembles (clean exponential decay), contrasting with the flat zero line on any associative quaternionic subalgebra $\HH \subset \OO$. (b)~Empirical quenched Lyapunov exponents $\hat{\lambda}_N$ ($\pm$ SE) stabilize cleanly well above the rigorous Minkowski/Flajolet--Sedgewick lower bounds $\log u_0^{(K)}$ of Theorem~\ref{thm:exp_decay}, illustrating the angular-cancellation gap highlighted in Open Problem~\ref{prob:lyap_limit}.}
\label{fig:lyapunov}
\end{figure}
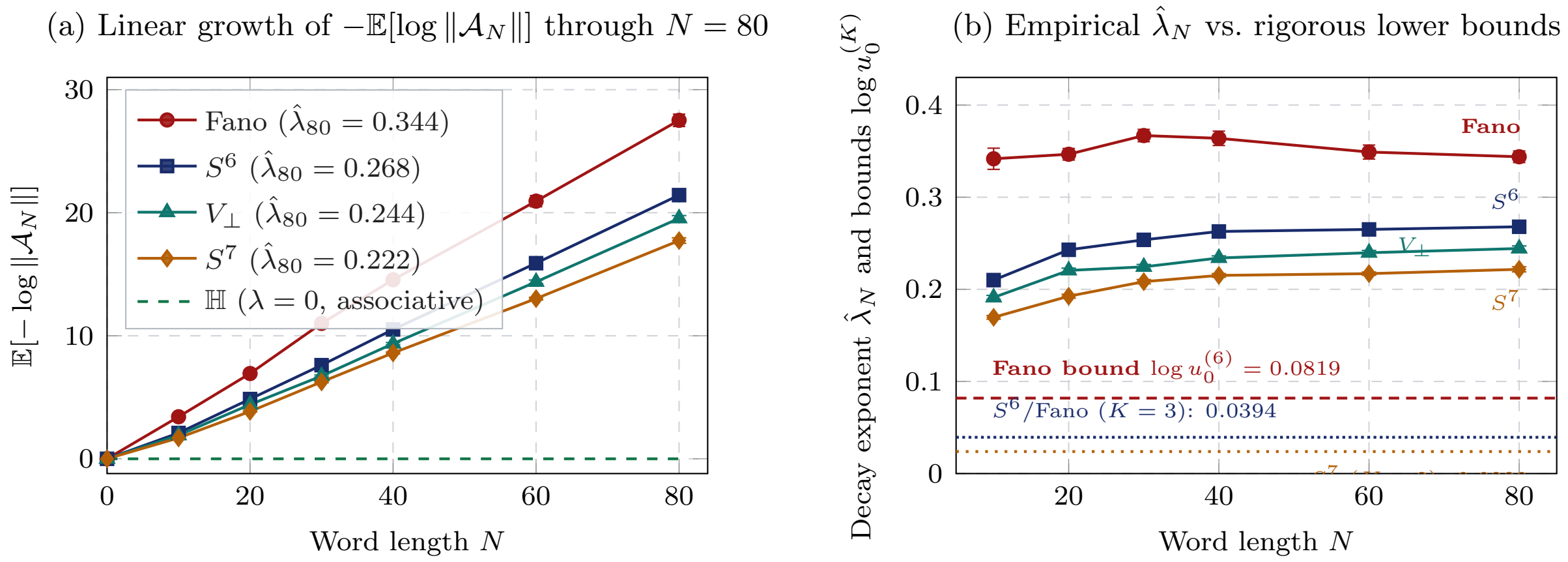

Three structural features stand out in Table~\ref{tab:lyapunov_num} and Figure~\ref{fig:lyapunov}:
\begin{enumerate}
    \item \textbf{Monotone Hierarchy Reflecting Associator Strength:} Across every $N \ge 10$, the quenched exponents obey the strict ordering
    \begin{equation}
    0 = \lambda_{\HH} \;<\; \lambda_{S^7} \approx 0.222 \;<\; \lambda_{V_\perp} \approx 0.244 \;<\; \lambda_{S^6} \approx 0.268 \;<\; \lambda_{\mathrm{Fano}} \approx 0.344.
    \end{equation}
    Including the associative real axis in $S^7$ reduces $\lambda$ relative to $S^6$; restricting to the $6$-dimensional transverse plane $V_\perp = e_7^\perp$ (associator Casimir $64$ vs.\ $96$ on $S^6$) lowers $\lambda$ correspondingly; and discrete Fano units maximize binary sign reversals ($\sigma_T \in \{\pm 1\}$) and exhibit the fastest decay.
    \item \textbf{Annealed vs.\ Quenched Gap:} Separate log-linear regressions on $\E[\|\calA_N\|^2]$ over $N \in [8, 32]$ yield clean exponential fits ($R^2 = 0.9992$ on $S^6$ and $R^2 = 0.9993$ on $S^7$, compared with $R^2 \approx 0.96$ for power laws), with annealed decay bases $\rho_{S^6} \approx 0.607$ and $\rho_{S^7} \approx 0.660$. Comparing with the quenched squares $e^{-2\hat{\lambda}_{80}} \approx 0.585$ ($S^6$) and $0.642$ ($S^7$) reveals a strict Jensen (large-deviation) gap $\rho > e^{-2\lambda}$, driven by rare nearly-quaternionic realizations whose norms decay more slowly than typical trajectories.
    \item \textbf{Complexified Split-Octonion Behavior:} Extending~\eqref{eq:dp_recurrence} $\C$-linearly to normalized Witt basis elements $A_a = \frac{1}{\sqrt{2}}(e_a + i e_{a'})$ and $A_a^* = \frac{1}{\sqrt{2}}(e_a - i e_{a'})$ in $\C \otimes \OO$ yields rapid exponential decay ($\E[\|\calA_N\|^2]$ falling from $0.465$ at $N = 3$ to $3.45 \times 10^{-5}$ at $N = 20$) for random Witt chains, whereas associative scalar chains $(1, \dots, 1)$ remain fixed at norm $1$.
\end{enumerate}

%==============================================================================
\section{Open Problems}
\label{sec:open}
%==============================================================================

Our results leave several sharp, self-contained mathematical problems open:

\begin{problem}[Existence and Exact Value of the Quenched Lyapunov Limit]
\label{prob:lyap_limit}
Theorem~\ref{thm:exp_decay} proves $\liminf_{N \to \infty} (-\frac{1}{N}\log\|\calA_N\|) \ge \log u_0^{(K)} > 0$ almost surely, and Table~\ref{tab:lyapunov_num} strongly indicates convergence to a deterministic limit $\lambda(\mu) = \lim_{N \to \infty} (-\frac{1}{N}\log\|\calA_N\|)$. Unlike associative random matrix products~\cite{furstenberg1960}, $\|\calA_N\|$ is not directly submultiplicative across arbitrary concatenations because Catalan trees couple left and right subwords across every split $k$. Proving almost-sure existence of $\lambda(\mu)$ (and determining whether $\lambda_{S^6}$ or $\lambda_{\mathrm{Fano}}$ admits a closed-form analytic characterization) remains open.
\end{problem}

\begin{problem}[Asymptotic Density of Exact Fano Zero-Cancellations]
By Corollary~\ref{cor:parity}, the fraction $p_0(N)$ of Fano basis words satisfying $\calA_N(w) = 0$ vanishes identically at every dyadic length $N = 2^k$, whereas $p_0(3) = 24/49 \approx 0.4898$, $p_0(5) = 3360/16807 \approx 0.1999$, and $p_0(6) = 8400/117649 \approx 0.0714$. What is the asymptotic decay rate of $p_0(N)$ along non-dyadic lengths $N \neq 2^k$?
\end{problem}

\begin{problem}[Non-Uniform Tree Measures and Higher Cayley--Dickson Algebras]
How does the non-associative Lyapunov exponent change when the uniform Catalan measure on $\calB_N$ is replaced by other natural tree ensembles (such as Galton--Watson or random binary search trees)? In the $16$-dimensional sedenions $\mathbb{S}$, where zero divisors break both~\eqref{eq:norm_mult} and Lemma~\ref{lem:assoc_kernel}, does the Universal Rigidity Theorem fail, and can non-associative averages grow or collapse faster than on $\OO$?
\end{problem}

\section*{Acknowledgments}
The author thanks colleagues for stimulating discussions on non-associative algebras and computational experimentation. In accordance with publisher transparency guidelines on generative AI, Google AI Studio was used as an interactive mathematical, computational, and editorial assistant during the development of this manuscript, including exploring the problem motivated by~\cite{smith2026octonion}, assisting with the formulation and cross-checking of proofs, implementing the exact $O(N^3)$ dynamic-programming and Monte Carlo computations, generating the vector figures, and drafting prose. The author reviewed, verified, and takes full scholarly responsibility for all theorems, proofs, numerical results, and conclusions.

\section*{Statements and Declarations}
\noindent\textbf{Funding:} No external funding was received for this work.\\
\noindent\textbf{Competing Interests:} The author declares no competing financial or non-financial interests.\\
\noindent\textbf{Data and Code Availability:} All exact enumerations and Monte Carlo benchmarks are reproducible directly from the $O(N^3)$ recurrence in Theorem~\ref{thm:bilinear_dp} and the Fano convention~\eqref{eq:fano_triads}.

\end{document}